\documentclass[twoside,11pt,openbib]{article}
\usepackage{latexsym}
\usepackage{amsmath}
\usepackage{amssymb}
\usepackage{amsfonts}
\usepackage{mathrsfs}

\numberwithin{equation}{section}
\newtheorem{thrm}{Theorem}[section]

\newtheorem{defn}{Definition}[section]
\newtheorem{lemm}{Lemma}[section]
\newtheorem{prop}{Proposition}[section]
\newtheorem{rmrk}{Remark}[section]

\newenvironment{prf}{{\em Proof:}}{\hfill{$\Box$}}

\newcommand{\alf}{\alpha}

\newcommand{\apx}{\approx}

\newcommand{\dert}{\frac{d}{dt}}

\newcommand{\Dlt}{{\Delta}}
\newcommand{\dt}[1]{\frac{d}{dt}{#1}}

\newcommand{\field}[1]{\mathbb{#1}}

\newcommand{\goto}{\rightarrow}
\newcommand{\gs}{\geqslant}

\newcommand{\hf}{\frac{1}{2}}

\newcommand{\lmd}{\lambda}

\newcommand{\lf}{\left}
\newcommand{\lra}[2]{\langle #1,#2\rangle}
\newcommand{\ls}{\leqslant}

\newcommand{\lsp}{\preccurlyeq}

\newcommand{\nbl}{\nabla}

\newcommand{\Omg}{\Omega}

\newcommand{\Pf}{\field P}

\newcommand{\prt}{\partial}

\newcommand{\R}{\field{R}}

\newcommand{\rt}{\right}

\newcommand{\T}{\field{T}}

\newcommand{\tht}{\theta}

\newcommand{\veps}{\varepsilon}

\begin{document}
\pagestyle{myheadings}

\title{\bf Global Regularity and Long-time Behavior of the Solutions to the 2D Boussinesq Equations without Diffusivity in a Bounded Domain}

\author{ Ning Ju\\
Department of Mathematics\\
401 Mathematical Sciences\\
Oklahoma State University\\
Stillwater, OK 74078, USA\\
Email: {\tt nju@okstate.edu}
}

\date{August 24th, 2015}

\maketitle

\begin{abstract}
 New results are obtained for global regularity and long-time behavior of the
 solutions to the 2D Boussinesq equations for the flow of an incompressible
 fluid with positive viscosity and zero diffusivity in a smooth bounded domain.
 Our first result for global boundedness of the solution $(u,\tht)$ in
 $D(A)\times H^1$ improves considerably the main result of the recent article
 \cite{Hu;Kukavica;Ziane:2013}. Our second result on global boundedness
 of the solution $(u,\tht)$ in $V\times H^1$ for both bounded domain and the
 whole space $\R^2$ is a new one. It has been open and also seems much more
 challenging than the first result. Global regularity of the solution
 $(u,\tht)$ in $D(A)\times H^2$ is also proved.\\

{\bf Keywords:} two dimensional dissipative Boussinesq equations, zero
 diffusivity, global regularity, long time behavior.

\end{abstract}

\markboth{Boussinesq Equations}{Ning Ju}

\indent
\baselineskip 0.6cm

\section{Introduction}
\label{s:intro}

\noindent

 Suppose that $\Omg\subset\R^2$ is a bounded domain with smooth boundary
 $\prt\Omg$. Consider the following Boussinesq equations for the flow of an
 incompressible fluid with positive viscosity and zero diffusivity:
\begin{equation}
 \label{e:u}
 \prt_t{u} - \Dlt u + u\cdot\nbl u + \nbl p = \tht e_2,\quad x\in\Omg,\ t>0
\end{equation}
\begin{equation}
 \label{e:c}
 \nbl\cdot u  =0,\quad x\in\Omg,\ t>0
\end{equation}
\begin{equation}
 \label{e:t}
\prt_t\tht + u\cdot\nbl\tht =0,\quad x\in\Omg,\ t>0,
\end{equation}
 where $u$ is the vector field of fluid velocity, $p$ and $\tht$ are the fluid
 pressure and density/temperature field and $e_2=(0,1)\in\R^2$. For simplicity
 of presentation, the viscosity constant has been set as $1$. The standard
 no-slip Dirichlet boundary condition is imposed:
\begin{equation}
\label{e:bc}
 u\big|_{\prt\Omg} = 0,\quad t>0.
\end{equation}
 The above equations are also equipped with initial conditions:
\begin{equation}
\label{e:ic}
 u(x,t)\big|_{t=0} = u_0(x),\quad \tht(x,t)\Big|_{t=0}=\tht_0(x),\quad x\in\Omg,
\end{equation}
 with proper $u_0$ and $\tht_0$ to be specified later.

 The 2D Boussinesq equations with {\em positive} viscosity and {\em positive}
 diffusivity is well known to be globally well-posed for both weak and strong
 solutions. This can be obtained following the classic theory for the 2D
 Navier-Stokes equations, cf. e.g. \cite{Constantin;Foias:1988} and
 \cite{Temam:1977}. See also a detailed proof of the global existence of the
 strong solutions as given in \cite{LiYG:2004}.
 If viscosity and diffusivity {\em both} vanish, then the global existence of
 strong solutions is still an open problem.

 Motivated by the works of \cite{Chae:2006} and \cite{Hou;Li:2005}, there have
 been extensive research activities on the global regularity of the solutions
 to the 2D Boussinesq equations for many cases sitting between the easiest one
 with full viscosity and diffusivity and the hardest one with zero viscosity
 and diffusivity. For simplicity of presentation, we refrain from discussing or
 citing extensive references for all those varying cases. In the following,
 we briefly discuss existing results just for the case to be studied in this
 article, i.e. the one with {\em positive and full} viscosity and {\em zero}
 diffusivity.
 The notations to be used in this article are basically standard and will be
 defined and explained in details in Section~\ref{s:pre}.

 First of all, existence and uniqueness of global weak solutions have been
 established in \cite{Danchin;Paicu:2008} for $\Omg=\R^2$ and in \cite{He:2012}
 for $\Omg$ being a smooth bounded domain. Global regularity for the solution
 $(u,\tht)\in H^1\times L^2$ can be easily proved. See e.g. \cite{Chae:2006},
 \cite{Danchin;Paicu:2011} and \cite{Hou;Li:2005} for $\Omg=\R^2$,
 \cite{Larios;Lunasin;Titi:2013} for $\Omg=\T^2$ and
 \cite{Hu;Kukavica;Ziane:2013} for $\Omg$ being a smooth bounded domain.
 See also Theorem~\ref{t:l2} and Theorem~\ref{t:v.l2}.

 For the case $\Omg=\R^2$, global in time regularity is proved in
 \cite{Chae:2006} for $(u,\tht)\in H^m\times H^m$, with the integer $m\gs 3$
 and in \cite{Hou;Li:2005} for $(u,\tht)\in H^m\times H^{m-1}$, with the integer
 $m\gs 3$. Global regularity is also proved in \cite{Danchin;Paicu:2008} for
 $(u,\tht)\in H^s\times H^{s-1}$, with real $s\gs 3$.

 For the case $\Omg$ being a bounded domain with a physical boundary, the
 particular difficulty in addition to complication of nonlinearity is the
 intensive vorticity along the boundary, which requires more careful treatment.
 Mathematically, it is convenient to use vorticity equation when $\Omg=\R^2$,
 as vortex stretching does not exist. Indeed, previous works for the case
 $\Omg=\R^2$ all take advantage of this feature. However, for a bounded domain
 $\Omg$, except the case of $\Omg=\T^2$, this approach is not as feasible since
 it is difficult to deal with vorticity along the boundary of $\Omg$.
 The case of $\Omg$ being a smooth bounded domain with no-slip boundary
 condition on $u$ has been studied in \cite{Lai;Pan;Zhao:2011}, the main result
 of which is the global regularity for solutions $(u,\tht)\in H^3\times H^3$
 under some {\em extra compatibility conditions} for initial datum.

 As pointed out by \cite{Hu;Kukavica;Ziane:2013}, due to the special feature of
 the equations, global in time regularity in the Sobolev spaces of intermediate
 orders is in fact much more difficult to be established than in the higher
 order Sobolve spaces $H^m\times H^m$ and $H^m\times H^{m-1}$ for $m\gs3$ as
 studied in the above mentioned literature. The main result of
 \cite{Hu;Kukavica;Ziane:2013} is global in time regularity of the solutions
 $(u,\tht)\in D(A)\times H^1$ for the case of a smooth bounded domain $\Omg$
 with no-slip boundary condition on $u$. Here, $A$ is the Stokes operator of
 $u$ with no-slip boundary condition on $\prt\Omg$. See Section~\ref{s:pre}
 for notations. This result is proved in \cite{Hu;Kukavica;Ziane:2013} using
 a new Gronwall lemma for two coupled nonlinear differential inequalities with
 one involving a logarithmic function. No additional compatibility conditions
 is required for $u_0$ beyond $u_0\in D(A)$. The analysis of
 \cite{Hu;Kukavica;Ziane:2013} also applies to the case of $\Omg=\R^2$ or $\T^2$
 with minor modifications to obtain similar results.
 It might be worthy of mentioning that for the case of $\Omg=\R^2$, global
 regularity in $H^2\times H^1$ is also obtained in \cite{Hou;Li:2005} using a
 boot-strapping argument, though much more complicated than the approach of
 \cite{Hu;Kukavica;Ziane:2013}. Notice also that vorticity equation is used in
 \cite{Hou;Li:2005}. Thus, it seems that the analysis of \cite{Hou;Li:2005} can
 not cover the case of a bounded domain with non-periodic boundary conditions.

 In this article, we study further the global regularity of solutions of this
 system and especially the long-time behavior of the solutions of this system.
 We will focus on the case of a {\em bounded} domain $\Omg\subset \R^2$. All
 the conclusions obtained in this article can be extended without adjustment to
 the case of $\Omg$ for which Poincar\'e inequality is valid. All the results
 except for those related to uniform upper bounds are also valid for the case
 $\Omg=\R^2$.
 
 First, a new proof of the $D(A)\times H^1$ global regularity will be presented
 using a method quite different from those of \cite{Hu;Kukavica;Ziane:2013}
 and \cite{Hou;Li:2005}. The new proof improves considerably the result of
 \cite{Hu;Kukavica;Ziane:2013} in the following sense: the best upper bounds of
 $\|Au(t)\|_2$ and $\|\nbl\tht(t)\|_2$ for $t\in[0,+\infty)$ that can be
 derived from the analysis of \cite{Hu;Kukavica;Ziane:2013} both grow
 {\em double exponentially} with respect to $t$ as $t\goto\infty$; while under
 the same conditions as those of \cite{Hu;Kukavica;Ziane:2013} and somewhat
 surprisingly, our analysis provides a {\em uniform} upper bound for
 $\|A(t)u\|_2$ with respect to $t\in[0,\infty)$ (and indeed also a {\em bounded
 absorbing set} for $u(t)$ in $D(A)$ as $t\goto\infty$) and an upper bound for
 $\|\nbl\tht(t)\|_2$ which is only {\em single exponential} with repsect to
 $t^2$ as $t\goto\infty$. See Theorem~\ref{t:u.h2} and Theorem~\ref{t:t.h1} for
 detailed statements of the result. Our new analysis seems also more
 transparent than previous ones. Notice that due to the inductive feature of
 the method of \cite{Hou;Li:2005}, it seems quite difficult to get upper bounds
 for $\|Au(t)\|_2$ and $\|\nbl\tht(t)\|_2$ in {\em explicit form of functions}
 of $t$ from the boot-strapping argument of \cite{Hou;Li:2005}. In this sense,
 the method of \cite{Hou;Li:2005} seems to be less sharp than that of
 \cite{Hu;Kukavica;Ziane:2013} and, as already mentioned above, it does not
 apply to the case of a bounded $\Omg$ with non-periodic boundary conditions.

 As our next main results, we prove the global regularity for $(u,\tht)\in V
 \times H^1$ and for $(u,\tht)\in D(A)\times H^2$. Here $V$ is the standard
 space of divergence-free vector fields in $L^2$ with their gradients also in 
 $L^2$. See Section~\ref{s:pre} for notations. Especially, the global regularity
 for $(u,\tht)\in V \times H^1$ seems to be quite unexpected complement of the
 existing results for global and local regularity. Notice that there are obvious
 essential reasons for previous analysis of regularity results to focus on the
 case $(u,\tht)\in H^2\times H^1$. For example, due to lack of diffusivity in
 the $\tht$ equation \eqref{e:t}, even though $\|\tht(t)\|_p$ is conserved for
 all $t\gs0$, there is no direct way to get global or even local in time
 regularity $\tht\in H^1$ just from $\tht$ equation \eqref{e:t} alone. Moreover,
 it seems still extremely difficult to obtain regularity of $\tht$ in $H^1$ by
 just coupling it with the enstropy equation of $u$, even just local in time.
 One key obstacle is that, in dealing with the involved nonlinearity,
 $\|u\|_{H^2}$ can not dominate $\|\nbl u\|_\infty$ in 2D domain via Sobolev
 imbedding. Even with one order higher regularity of $u$ coupled, it is still
 technically complicated to combine the tool of Brezis-Gallouet type
 inequalities to deal with regularity of $\tht$ in $H^1$, as can be seen in
 \cite{Hu;Kukavica;Ziane:2013}, \cite{Hou;Li:2005} and in the proof of the
 first main result of this article as given in Section~\ref{s:h2.h1}. Now, with
 only $H^1$ regularity for $u$ coupling with $H^1$ regularity for $\tht$,
 proving global regularity of $(u,\tht)\in V\times H^1$ appears {\em much} more
 challenging than proving global regularity of $(u,\tht)\in D(A)\times H^1$,
 since Brezis-Gallouet type inequalities would {\em not} help any more.
 Indeed, to the knowledge of the author, existence of $V\times H^1$ solutions,
 has been {\em open} for any kind of domains, even though the essential
 difficulty here is indeed only a {\em local} property. In this article, some
 new ideas from spectral decomposition analysis will be used to resolve
 completely the essential difficulty. The main result on global regularity of
 $(u,\tht)\in V\times H^1$ thus obtained will be summarized in
 Theorem~\ref{t:h1}. As the second application of spectral decomposition and
 for completeness, we prove our next main result, Theorem~\ref{t:h2}, for the
 global regularity of $(u,\tht)\in D(A)\times H^2$. Since we will not use
 Brezis-Gallouet type inequalities here, the proof is a lot simpler than
 otherwise.

 We comment that our proofs of the above mentioned global regularity results are
 self-contained and global regularity results in higher order Sobolev spaces
 can thus be obtained as consequence of our global regularity results.

 The rest of this article is organized as following:\\
 In Section~\ref{s:pre}, we give the notations, briefly review the background
 results and recall some important facts crucial to later analysis.
 In Section~\ref{s:h2.h1}, we state and prove Theorem~\ref{t:u.h2} and
 Theorem~\ref{t:t.h1} for the global regularity of the solutions
 $(u, \tht)\in D(A)\times H$.
 In Section~\ref{s:h1.h1}, we state and prove Theorem~\ref{t:h1} for the global
 regularity of the solutions  $(u, \tht)\in V\times H^1$.
 In Section~\ref{s:h2.h2}, we state and prove Theorem~\ref{t:h2} for the global
 regularity of the solutions  $(u, \tht)\in D(A)\times H^2$.

\section{Preliminaries}
\label{s:pre}

\subsection{Notations and Some Basic Results}

 Throughout this article, we use the notations that, for real numbers $A$ and
 $B$,
\[  A \lsp B \quad \text{ iff }\quad  A\ls C\cdot B,\]
 and
\[  A \apx B \quad \text{ iff }\quad  c\cdot A\ls B \ls C\cdot A, \]
 for some positive constants $c$ and $C$ independent of $A$ and $B$.

 Recall that $\Omg$ is a bounded domain in $\R^2$ with sufficiently smooth
 boundary $\prt\Omg$. By a domain, we always mean an open and connected open
 subset of $\R^2$. Denote by $L^p(\Omg)$ ($1\ls p <+\infty$) the classic
 Lebesgue $L^p$ spaces with the norms:
\[ \|f\|_p := \left( \int_\Omg|f(x)|^p\ dx \right)^\frac{1}{p},\quad
  \forall f\in L^p(\Omg). \]
 The Lebsgue space $L^\infty(\Omg)$ is defined as the Banach space of (classes
 of) real functions on $\Omg$ which are measurable and essentially bounded with
 the norm
\[ \|f\|_{\infty} := \text{ess\,}\sup_{x\in\Omg}|f(x)|. \]
 Denote by $H^m(\Omg)$ ($m\gs1$) the Sobolev space for square-integrable
 functions with square-integrable weak derivatives up to order $m$ with the norm
\[ \|f\|_{H^m} := \lf(\sum_{|\alf|\ls m}\|D^\alf f\|_2^2\rt)^\hf,\quad
 \forall f\in H^m(\Omg). \]
 We use the standard notations of the following functional spaces for the
 initial value problem with no-slip boundary condition:
\[ H:= \lf\{ v\in (L^2(\Omg))^2\ \big| \  \nbl\cdot v =0\ \text{ in }\Omg,\quad
 v\cdot n = 0 \ \text{ on } \prt\Omg, \rt\}, \]
\[ V:= \lf\{ v\in (H_0^1(\Omg))^2\ \big|\ \nbl\cdot v = 0\ \text{ in }\Omg.
   \rt\}, \]
 where the equalities are in distribution sense, $n$ is the unit vector
 normal to $\prt\Omg$ and outward with respect to $\Omg$ and $H_0^1(\Omg)$ is
 the closure of $C_0^\infty(\Omg)$ in $H^1(\Omg)$, i.e. the space of the
 functions in $H^1$ with zero-trace on $\prt\Omg$. Later on, we do not
 distinguish the notations for vector and scalar function spaces, which are
 self-evident from the context. Therefore, we will use e.g. $L^p(\Omg)$ and
 $H^m(\Omg)$ or simply $L^p$ and $H^m$ to denote $(L^p(\Omg))^2$ and
 $(H^m(\Omg))^2$ respectively provided there is no confusion.
 
 Denote by $\Pf$ the Helmholtz orthogonal projector in $(L^2(\Omg))^2$ onto
 $H$. The Stoke operator $A: D(A)\subset H \mapsto H$ is defined as
\[ A:=-\Pf\Dlt,\quad D(A)= (H^2(\Omg))^2\cap V.\]
 Recall that
\begin{equation}
\label{e:lp}
 \|\Pf v\|_H \lsp \|v\|_2,\ \forall v\in L^2(\Omg)^2, \quad
   \|\nbl\Pf v\|_2 \lsp \|v\|_{H^1},\ \forall v\in H^1(\Omg)^2.
\end{equation}
 By the spectral theorem for $A$, there exists a sequence $\{\lmd_j\}_1^\infty$
\[ 0<\lmd_1\ls \cdots \ls \lmd_j \ls \lmd_{j+1} \ls \cdots,\quad
   \lim_{j\goto\infty}\lmd_j = \infty, \]
 and a family $\{ w_j\}_1^\infty\subset D(A)$ which is orthonormal in $H$ such
 that
\[ A w_j = \lmd_j w_j,\quad \forall j =1, 2, \dots. \]
 Denote, for $\alf\in[0,\infty)$,
\[ A^\alf v := \sum_{j=1}^\infty \lmd_j^\alf(v,w_j)w_j,\quad
  \forall v\in D(A^\alf)(A), \]
 where $(\cdot,\cdot)$ is the inner product in $H$ and
\[ D(A^\alf) := \lf\{ v\in H \ \Big| \
   \sum_{j=1}^\infty \lmd_j^{2\alf}|(v,w_j)|^2<\infty \rt\}.\]
 Notice that, $D(A^\hf)=V$ and
\begin{equation}
\label{e:a.hf}
 \|A^\hf v\|_2 =\|\nbl v\|_2,\quad \forall v\in V.
\end{equation}
 Moreover, $\|A\cdot\|_2$ is equivalent to $\|\cdot\|_{H^2}$ on $D(A)$, i.e.
\[   \|A v\|_2 \apx \| v \|_{H^2} \quad \forall v\in D(A). \]
 For $0\ls r\ls s \ls 1$, there is a Poincar\'e inequality:
\[  \lmd_1^{s-r}\|A^r u\|_2 \ls \|A^s u\|_2. \]
 Define the family of spectral projection operators $\{P_m\}_1^\infty$ as
\[ P_1 v := \sum_{\lmd_j<2} (u, w_j)w_j, \quad
   P_m v := \sum_{2^{m-1}\ls \lmd_j <2^m} (u,w_j)w_j,\quad
    \text{ for }\quad m\gs 2. \]
 Then,
\begin{equation}
\label{e:a.eqv}
   \|P_mA^s v\|_2 = \|A^s P_m v\|_2 \apx 2^{ms}\|P_m v\|_2,
  \quad \text{ for } m\gs 1,\ s\gs 0.
\end{equation}
 Introducing the bi-linear operator,
\[ B(u,v) := \Pf(u\cdot\nbl v),\quad \forall u, v\in V, \]
 which maps $V\times V$ into $V'$ (the dual space of $V$), one can also
 reformulate equations \eqref{e:u}-\eqref{e:c} as
\begin{equation}
\label{e:u.v}
 u_t + Au + B(u,u) =\Pf(\tht e_2).
\end{equation}
 It can be shown that (see e.g. \cite{Hu;Kukavica;Ziane:2013})
\begin{equation}
\label{e:b}
 \|\nbl B(u,v)\|_2
 \lsp \|u\|_2^{\frac{1}{4}}\|Au\|_2^{\frac{3}{4}}
      \|v\|_2^{\frac{1}{4}}\|Av\|_2^{\frac{3}{4}}
   + \|u\|_2^\hf\|Au\|_2^\hf\|Au\|_2.
\end{equation}
  
\begin{defn}
 Given $(u_0,\tht_0)\in H\times L^2$, $(u(t),\tht(t))$ is a {\em (global) weak
 solution} of the initial boundary value problem \eqref{e:u}-\eqref{e:ic} if
 for any $T>0$, the following statements are valid:
\begin{enumerate}
 \item $\tht \in C([0,T]; H)$, $v\in C([0,T]) \cap L^2(0, T; V)$.
 \item For any $\psi\in(C^1([0,T]\times\Omg))^2$ such that $\nbl\cdot\psi=0$,
\[ \int_0^T \Big[(u, \prt_t\psi)-(\nbl u,\nbl\psi)
 -(u\cdot\nbl u-\tht e_2, \psi)\Big]\ dt  = (u(T),\psi(T))-(u_0,\psi(0)).\]
 \item For any $\phi\in C^1([0,T]\times\Omg)$,
\[\int_0^T\Big[(\tht,\prt_t\phi) - (u\cdot\phi)\Big]\ dt
  = (\tht(T),\phi(T))-(\tht_0,\phi(0)). \]
\end{enumerate}
\end{defn}

 Notice that one could define a weak solution even somewhat {\em weaker} than
 the above one. We adopt this one since existence and especially uniqueness of
 the weak solution to the problem \eqref{e:u}-\eqref{e:ic} has been proved in
 \cite{Danchin;Paicu:2008} for the case when $\Omg=\R^2$, and in \cite{He:2012}
 for the case of a smooth ($C^2$) bounded domain.

 As mentioned in Section~\ref{s:intro}, for $\Omg=\R^2$, global regularity of
 the solutions in $H^m\times H^{m}$ with integer $m\gs3$ is proved in
 \cite{Chae:2006} and global regularity in $H^m\times H^{m-1}$ with integer
  $m\gs3$ is proved in \cite{Hou;Li:2005}, which is extended in
 \cite{Danchin;Paicu:2008} for $(u,\tht)\in H^s\times H^{s-1}$, with real
 $s\gs 3$. In \cite{Lai;Pan;Zhao:2011} global regularity in $H^3\times H^3$ is
 proved when $\Omg$ is a smooth {\em bounded} domain, under {\em extra}
 compatibility conditions for initial data. In the more recent work
 \cite{Hu;Kukavica;Ziane:2013}, global regularity in $D(A)\times H^1$ is
 obtained where $\Omg$ can be either a bounded domain, $\T^2$ or $\R^2$. No
 {\em extra} compatibility conditions for initial data other than $u_0\in D(A)$
 is needed.
 
 The main goal of this article is to study global regularity and long-time
 behavior of the solutions for the case of a bounded domain. The main results
 of this article, as discussed in Section~\ref{s:intro}, are
 Theorem~\ref{t:u.h2}, Theorem~\ref{t:t.h1}, Theorem~\ref{t:h1} and
 Theorem~\ref{t:h2}. The first two improve considerably the main result of
 \cite{Hu;Kukavica;Ziane:2013}. The last two are new results which seem to
 have been open.

 For convenience of reading, we recall the following formulation of the
 Uniform Growall Lemma, which will be used frequently in later discussion.
 The detailed proof of the lemma can be found e.g. in \cite{Temam:1988}.
\begin{lemm}[Uniform Gronwall Lemma]
\label{l:ug}
 Let $g$, $h$ and $y$ be three non-negative locally integrable functions on
 $(t_0, +\infty)$ such that
\[ \frac{dy}{dt} \ls gy + h, \qquad \forall t\gs t_0, \]
 and
\[ \int_t^{t+r}g(s)\ ds \ls a_1, \qquad \int_t^{t+r}h(s)\ ds \ls a_2, \qquad 
   \int_t^{t+r}y(s)\ ds \ls a_3, \qquad \forall t\gs t_0, \]
 where $r$, $a_1$, $a_2$ and $a_3$ are positive constants. Then
\[ y(t+r) \ls \left(\frac{a_3}{r}+a_2\right)e^{a_1}, \qquad \forall t\gs t_0. \]
\end{lemm}

\subsection{Global regularity in $H\times L^2$}
\label{ss:h.l2}

 For convenience of later discussion, we recall the following global regularity
 result:
\begin{thrm}
\label{t:l2}
 Suppose $(u_0,\tht_0)\in H\times L^2$. Then, \eqref{e:u}-\eqref{e:ic} has a
 unique weak solution $(u,\tht)$. Moreover, the following are valid:
\begin{equation}
\label{e:t.l2}
\|\tht(t)\|_2=\|\tht_0\|_2,\quad \forall t>0,
\end{equation}
\begin{equation}
\label{e:u.l2}
 \|u(t)\|_2^2 \ls \|u_0\|_2^2e^{-\lmd_1 t}
  + \lf(\frac{\|\tht_0\|_2}{\lmd_1}\rt)^2\lf(1-e^{-\lmd_1 t}\rt),\quad
 \forall t>0.
\end{equation}
\begin{equation}
\label{e:u.nbl.avg}
\int_t^{t+r}\|A^\hf u(\tau)\|_2^2\ d\tau
 \ls \|u_0\|_2^2e^{-\lmd_1 t}
  + \lf(\frac{\|\tht_0\|_2}{\lmd_1}\rt)^2\lf(1-e^{-\lmd_1 t}\rt)
  + \frac{r}{\lmd_1}\|\tht_0\|_2^2,\quad \forall t,r>0.
\end{equation}
\end{thrm}

\begin{prf}

 Existence and uniqueness of the weak solution have already been established in
 \cite{He:2012}.
 For completeness, we provide the brief formal {\em a priori} estimates yielding
 \eqref{e:t.l2}-\eqref{e:u.nbl.avg}. Due to existence and uniqueness of the weak
 solution, these a priori estimates can be justified rigorously via standard
 argument.
 First of all, \eqref{e:t.l2} is an immediate consequence of \eqref{e:t},
 \eqref{e:c} and \eqref{e:bc} via the following simple computations:
\begin{equation*}
\begin{split}
 \frac{1}{2}\dert{\|\tht\|_2^2}
&=\int_\Omg \tht_t\tht\ dx = - \int_\Omg (u\cdot\nbl\tht)\tht\ dx
 =-\frac{1}{2}\int_\Omg u\cdot\nbl (\tht^2)\ dx\\
&= \frac{1}{2}\int_{\prt\Omg}(u\cdot n)\tht^2\ dS
 -\frac{1}{2}\int_\Omg(\nbl\cdot u)\tht^2\ dx = 0.
\end{split}
\end{equation*}
 Next, taking inner product of \eqref{e:u} with $u$ and using \eqref{e:c} and
 \eqref{e:bc}, we obtain
\[ \hf\dert{\|u\|_2^2} + \|A^\hf u\|_2^2 \ls\|\tht\|_2\|u\|
 \ls \frac{\lmd_1}{2}\|u\|_2^2 + \frac{1}{2\lmd_1}\|\tht\|_2^2
 \ls \hf\|A^\hf u\|_2^2 + \frac{1}{2\lmd_1}\|\tht\|_2^2. \]
 Therefore,
\begin{equation}
\label{e:u.l2.nbl}
\dert{\|u\|_2^2} + \|A^\hf u\|_2^2 \ls \frac{1}{\lmd_1}\|\tht_0\|_2^2,
\end{equation}
 implying that
\[ \dert{\|u\|_2^2} +\lmd_1\|u\|_2^2 \ls \frac{1}{\lmd_1}\|\tht_0\|_2^2, \]
 from which \eqref{e:u.l2} follows by a direct integration. Now,
 \eqref{e:u.nbl.avg} can be obtained from \eqref{e:u.l2.nbl} and \eqref{e:u.l2}.

\end{prf}

\begin{rmrk}
\label{r:abs.H}
 We see from \eqref{e:u.l2} that, for any fixed $\|\tht_0\|_2$, there is a
 {\em bounded} absorbing set for $u$ in $H$ with the radius being
 $\|\tht_0\|_2/\lmd_1$ and not depending on $\|u_0\|_2$. It is also easy to
 see from \eqref{e:u.nbl.avg} that there is a {\em bounded} absorbing set of
 $\int_t^{t+r}\|A^\hf u\|_2^2d\tau$ in $\R_+$ with fixed $\|\tht_0\|_2$ and
 $r(>0)$.
\end{rmrk}

\subsection{Global Regularity in $V\times L^2$}
\label{ss:u.V}

 The following global regularity will also be used in later discussion:
\begin{thrm}
\label{t:v.l2}
 Suppose that $(u_0,\tht_0)\in V\times L^2$. Then, there is a
 $t_0=t_0(u_0,\tht_0)>0$, such that the unique weak solution
 $(u,\tht)$ of \eqref{e:u}-\eqref{e:ic} satisfies
\begin{equation}
\label{e:u.nbl.loc}
 \|A^\hf u(t)\|_2^2 \ls 1 + 2\|A^\hf u_0\|_2^2,\quad \forall t\in[0, t_0].
\end{equation}
 Moreover, for $t\gs t_0$, the following uniform bound of $\|A^\hf u(t)\|_2$
 is valid:
\begin{equation}
\label{e:u.nbl.glb}
 \|A^\hf u(t)\|_2^2 \ls \lf[\|u_0\|_2^2e^{-\lmd_1t} t_0^{-1}
 +(\lmd_1^{-1}+t_0)\|\tht_0\|_2^2
 + \frac{\|\tht_0\|_2^2}{\lmd_1^2t_0}\lf(1-e^{-\lmd_1t}\rt)\rt] e^{a_1},
\end{equation}
 where, there is a generic constant $c_0>0$, such that
\begin{equation}
\label{e:a1}
 a_1 \ls c_0\lf( \|u_0\|_2^2e^{-\lmd_1t} + \frac{\|\tht_0\|_2^2}{\lmd_1^2}\rt)
 \lf[\|u_0\|_2^2e^{-\lmd_1t}+\frac{\|\tht_0\|_2^2}{\lmd_1^2}\lf(1-e^{-\lmd_1t}\rt)
 + \frac{t_0}{\lmd_1}\|\tht_0\|_2^2\rt].
\end{equation}
 Furthermore, for any $t,r>0$,
\begin{equation}
\label{e:u.h2.avg}
\begin{split}
 \int_t^{t+r}\|Au(\tau)\|_2^2\ d\tau
 \lsp& \|A^\hf u(t)\|_2^2 + r\|\tht_0\|_2^2\\
 &+ \int_t^{t+r}\lf[\|u_0\|_2^2e^{-\lmd_1\tau}
 +\lf(\frac{\|\tht_0\|_2}{\lmd_1}\rt)^2\rt]\|A^\hf u(\tau)\|_2^4\ d\tau.
\end{split}
\end{equation}
\end{thrm}

\begin{prf}
 Due to Theorem~\ref{t:l2}, the following {\em a priori} estimates can all be
 justified rigorously. Taking the inner product of \eqref{e:u.v} with $Au$, we
 obtain
\begin{equation}
\label{e:u.nbl.l2}
\begin{split}
 \hf\dert{\|A^\hf u\|_2^2} + \|Au\|_2^2
&\ls \|\tht\|_2\|Au\|_2 +|\lra{B(u,u)}{Au}|\\
&\ls \|\tht_0\|_2\|Au\|_2 + C\|u\cdot\nbl u\|_2\|Au\|_2\\
&\ls \|\tht_0\|_2\|Au\|_2
 + C\|u\|_2^\hf\|\nbl u\|\|Au\|_2^\frac{3}{2}\\
&\ls \hf\|Au\|_2^2
 + C\lf(\|\tht_0\|_2^2 + \|u\|_2^2\|A^\hf u\|_2^4\rt),
\end{split}
\end{equation}
 where we have used Gagliardo-Nirenberg inequality. Therefore,
 by \eqref{e:u.nbl.l2} and \eqref{e:u.l2}, we have
\begin{equation*}
\begin{split}
 \dert{\|A^\hf u\|_2^2} + \|Au\|_2^2 
&\lsp \|\tht_0\|_2^2
 + \lf[\|u_0\|_2^2 +\lf(\frac{\|\tht_0\|_2}{\lmd_1}\rt)^2\rt]\|A^\hf u\|_2^4.
\end{split}
\end{equation*}
 Therefore,
\begin{equation*}
\dert{\|A^\hf u\|_2^2} \ls C_0 + C_1\|A^\hf u\|_2^4,
\end{equation*}
 where
\[ C_0=C\|\tht_0\|_2^2,\quad
 C_1=C\lf[\|u_0\|_2^2 +\lf(\frac{\|\tht_0\|_2}{\lmd_1}\rt)^2\rt]. \]
 Let $y(t)=1+\|A^\hf u(t)\|_2^2$ and $C_2=\max\{C_0,C_1\}$. Then,
 $y'(t) \ls  C_2 y^2$. Integrating the inequality with respect to $t$ yields
\[ y(t) \ls \frac{y(0)}{\lf(1-C_2y(0)t\rt)}.\]
 Choose $t_0>0$ such that
\[ t_0 \ls \frac{1}{2C_2y(0)} = \frac{1}{2C_2\lf(1+\|A^\hf u_0\|_2^2\rt)}. \]
 Then, for $t\in[0, t_0]$,
\begin{equation*}
 1 + \|A^\hf u(t)\|_2^2
 \ls \frac{1+\|A^\hf u_0\|_2^2}{1-C_2\lf(1+\|A^\hf u_0\|_2^2\rt)t}
 \ls 2\lf(1+\|A^\hf u_0\|_2^2\rt).
\end{equation*}
 This proves the local estimate \eqref{e:u.nbl.loc}. Now we prove the global
 estimate \eqref{e:u.nbl.glb}. By \eqref{e:u.nbl.l2} and \eqref{e:u.l2}, we
 also have
\begin{equation}
\label{e:u.nbl.l2.1}
\begin{split}
 \dert{\|A^\hf u\|_2^2} + \|Au\|_2^2 
&\lsp \|\tht_0\|_2^2
 + \lf[\|u_0\|_2^2e^{-\lmd_1t}
 +\lf(\frac{\|\tht_0\|_2}{\lmd_1}\rt)^2\rt]\|A^\hf u\|_2^4,
\end{split}
\end{equation}
 Noticing \eqref{e:u.nbl.avg}, we can apply the Lemma~\ref{l:ug} to
 \eqref{e:u.nbl.l2.1} with $r=t_0$ to obtain the global estimate
 \eqref{e:u.nbl.glb}-\eqref{e:a1}. Finally, \eqref{e:u.h2.avg} follows from
 integrating \eqref{e:u.nbl.l2.1} with respect to $t$ from $t$ to $t+r$.

\end{prf}

\begin{rmrk}
\label{r:abs.V}
 We see from \eqref{e:u.nbl.glb} and \eqref{e:a1} that, for any fixed
 $\|\tht_0\|_2$, there is a {\em bounded} absorbing set for $u$ in space $V$
 with the radius depending on $\|\tht_0\|_2$ but not $u_0$. It is also easy to
 see from \eqref{e:u.h2.avg} that there is a {\em bounded} absorbing set of
 $\int_t^{t+r}\|Au\|_2^2d\tau$ in $\R_+$ for fixed $\|\tht_0\|_2$ and $r(>0)$
 which is also independent of $u_0$.
\end{rmrk}

\section{Global Regularity in $D(A)\times H^1$}
\label{s:h2.h1}

 The first main result of this article is achieved via an alternative proof of
 the global regularity of the solution $(u,\tht) \in D(A)\times H^1$. Proving
 this global regularity is the main result of the recent work
 \cite{Hu;Kukavica;Ziane:2013}. The improvement of our new result over that of
 \cite{Hu;Kukavica;Ziane:2013}\footnote{The approach of \cite{Hou;Li:2005}
 applies only to the case of $\Omg$ with no boundary and yields less sharp
 estimate than that of \cite{Hu;Kukavica;Ziane:2013}.} is in the following
 sense: the {\em best} upper bounds for $\|Au(t)\|_2$ and $\|\nbl\tht(t)\|_2$
 that can be obtained from the proof of \cite{Hu;Kukavica;Ziane:2013} both
 increase {\em double exponentially} with respect to $t$ as $t\goto\infty$;
 while under the same conditions, our new proof provides a {\em uniform} upper
 bound for $\|Au(t)\|_2$ with respect to $t\in[0,\infty)$ and an upper bound for
 $\|\nbl\tht(t)\|_2$ which grows only {\em single exponentially} with respect
 to $t^2$. Indeed, we also obtain a bounded absorbing set for $\|Au(t)\|_2$ in
 $\R_+$, which depends only on $\tht_0$ but not on $u_0$. To achieve the
 improvement, our method has to be quite different from the previous ones. We
 split the statement of our first main result into Theorem~\ref{t:u.h2} and
 Theorem~\ref{t:t.h1}. We also split the proof of our first main result into
 the following three subsections.

\subsection{Local Regularity in $D(A)\times H^1$}

 As the preparation for the analysis in the next two subsections, we prove the
 following local in time boundedness of $\|Au(t)\|_2$ and $\|\nbl\tht(t)\|_2$:
\begin{prop}
\label{p:h2.h1.loc}
 Suppose $(u_0,\tht_0)\in D(A)\times H^1$. Then, there exists a
 $t_1\in(0, t_0]$, such that
\begin{equation}
\label{e:h2.h1.loc}
 \|Au\|_2^2 + \|\nbl\tht\|_2^2
 \ls 2\lf(\|Au_0\|_2^2 + \|\nbl\tht_0\|_2^2\rt), \quad \forall t\in[0,t_1],
\end{equation}
 where $t_0>0$ is given in \eqref{e:u.nbl.loc}.
\end{prop}

\begin{prf}

 Taking inner product of \eqref{e:u.v} with $A^2u$ and using \eqref{e:b} and
 \eqref{e:lp}, we have
\begin{equation}
\label{e:u.h2}
\begin{split}
 \dert{\|Au\|_2^2} + \|A^\frac{3}{2}u\|_2^2
&\lsp \|\tht\|_{H^1}^2 + \|u\|_2\|Au\|_2^3\\
&\lsp \|\tht_0\|_2^2 + \|\nbl\tht\|_2^2+ \|u\|_2\|Au\|_2^3.
\end{split}
\end{equation}
 Applying $\nbl$ to \eqref{e:t} and taking the inner product of the derived
 equation with $\nbl$, we have
\begin{equation}
\label{e:t.nbl}
 \hf\dert{\|\nbl\tht\|_2^2}
 = -\int_\Omg \sum_{i,j=1}^2\prt_j\tht\prt_ju_i\prt_i\tht\ dx
 \ls \|\nbl u\|_\infty\|\nbl\tht\|_2^2,
\end{equation}
 where we have used an integration by parts, no-slip boundary condition
 and the fact that $\nbl\cdot u=0$. Therefore,
\begin{equation}
\label{e:t.h1}
 \dert{\|\nbl\tht\|_2^2} 
\ls \hf\|A^\frac{3}{2} u\|_2^2 + C\|\nbl\tht\|_2^4.
\end{equation}
 Adding \eqref{e:u.h2} and \eqref{e:t.h1} yields
\[ \dert{\lf(\|Au\|_2^2 + \|\nbl\tht\|_2^2\rt)}
  + \|A^\frac{3}{2} u\|_2^2 \ls C(1+\|\tht_0\|_2^2)
 +C\lf(\|Au\|_2^2 + \|\nbl\tht\|_2^2\rt)^2. \]
Therefore, there exists a $t_1\in(0, t_0]$ such that \eqref{e:h2.h1.loc} is
 valid.  This finishes the proof of local regularity of
 $(u,\tht)\in L^\infty(0, t_1; D(A)\times H^1)$.

\end{prf}

\subsection{Global Uniform Estimates of $\|A u\|_2$ and $\|u_t(t)\|_2$}

 In this subsection, we prove the global in time uniform boundedness of
 $\|Au(t)\|_2$. Our approach is completely different from previous methods used
 in \cite{Hu;Kukavica;Ziane:2013}, which can only get an upper bound of
 $\|Au(t)\|_2$ which gowths double exponentially with respect to $t$ as
 $t\goto\infty$. 

\begin{thrm}
\label{t:u.h2}
 Suppose $(u_0,\tht_0)\in(D(A), H^1)$. Then,
 \[ u\in L^\infty(0, +\infty; D(A)),\quad u_t\in L^\infty(0, +\infty; H)\cap
 L^2_{loc}(0, +\infty; V). \]
 Moreover, for fixed $\|\tht_0\|_{H^1}$, there exists a bounded absorbing set
 in $\R_+$ for $\|Au(t)\|_2$, $\|u_t(t)\|_2$ and
 $\int_t^{t+1}\|u_t(\tau)\|_V^2 d\tau $, which is independent of $u_0$.
\end{thrm}

\begin{prf}

 To prove uniform boundedness of $u\in L^\infty(0,\infty;D(A))$, it is enough
 to prove uniform boundedness of $u_t\in L^\infty(t_1,\infty;H)$ for
 sufficiently small $t_1>0$. However, we will prove in the following slightly
 stronger uniform boundedness of $u_t\in L^\infty(0,\infty;H)$, {\em without}
 imposing the condition that $u_t(0)\in H$.

{\bf Step 1}. Local regularity of  $u_t\in L^\infty(0, t_1;H)$.

 By \eqref{e:u.v}, we have
\begin{equation*}
\|u_t\|_2^2
 = \lra{- Au - B(u,u) +{\field P}\tht e_2}{u_t}
 \ls (\|A u\|_2 + \|u\cdot\nbl u\|_2 + \|\tht\|_2)\|u_t\|_2.
\end{equation*}
 Therefore,
\begin{equation}
\label{e:u.t.l2}
 \|u_t\|_2 \lsp \|Au\|_2(1 +\|\nbl u\|_2) + \|\tht_0\|_2.
\end{equation}
 Hence, by \eqref{e:u.t.l2}, \eqref{e:u.nbl.loc} and \eqref{e:h2.h1.loc},
 we have
\begin{equation*}
\|u_t\|_2 
\lsp (\|Au_0\|_2+\|\nbl\tht_0\|)(1+\|\nbl u_0\|_2) + \|\tht_0\|_2,\quad
 \forall t\in(0,t_1],
\end{equation*}
 which proves local regularity of $u_t\in L^\infty(0, t_1;H)$.

{\bf Step 2}. Global regularity of $u_t\in L^\infty(0,\infty; H)\cap
 L^2(0,\infty;V)$.

 First, square \eqref{e:u.t.l2} and integrate with respect to $t$ to get
\begin{equation}
\label{e:u.t.avg}
 \int_t^{t+r}\|u_t(\tau)\|_2^2\ d\tau \lsp  r\|\tht_0\|_2^2
 + \int_t^{t+r} \|Au(\tau)\|_2^2(1+\|\nbl u(\tau)\|_2^2)\ d\tau.
\end{equation}
 By \eqref{e:u.t.avg}, Theorem~\ref{t:l2} and Theorem~\ref{t:v.l2}, we see that
 $\int_t^{t+r}\|u_t(\tau)\|_2^2\ d\tau$ is uniformly bounded.
 By Remark~\ref{r:abs.H} and Remark~\ref{r:abs.V}, it is also easy to see that
 there is a bounded absorbing set in $\R_+$ for
 $\int_t^{t+r}\|u_t(\tau)\|_2^2\ d\tau$ when $\tht_0$ and $r>0$ are fixed.
 Next, apply $\prt_t$ to \eqref{e:u} to get
\begin{equation}
\label{e:u.t}
 \prt_t u_t -\Dlt u_t +u_t\cdot\nbl u +u\cdot\nbl u_t +\nbl p_t = \tht_t e_2.
\end{equation}
 Taking inner product of \eqref{e:u.t} with respect $u_t$ yields
\begin{equation*}
\begin{split}
 \hf\dert{\|u_t\|_2^2} + \|\nbl u_t\|_2^2 &=-\lra{u_t\cdot\nbl u}{u_t}
 -\lra{u\cdot\nbl u_t}{u_t}-\lra{p_t}{u_t} +\lra{\tht_t e_2}{u_t}\\
&=-\lra{u_t\cdot\nbl u}{u_t} -\lra{u\cdot\nbl \tht e_2}{u_t}\\
&=-\lra{u_t\cdot\nbl u}{u_t}+ \lra{\tht u}{\nbl u_{2,t}}\\
&\ls \|\nbl u\|_2\| u_t\|_4^2 + \|u\tht\|_2^2 + \frac{1}{4}\|\nbl u_t\|_2^2\\
&\ls C\|\nbl u\|_2^2\|u_t\|_2^2 + C\|Au\|_2^2\|\tht_0\|_2^2
 + \hf\|\nbl u_t\|_2^2.
\end{split}
\end{equation*}
 Therefore,
\begin{equation}
\label{e:u.t.2}
 \dert{\|u_t\|_2^2} + \|\nbl u_t\|_2^2
 \lsp \|\nbl u\|_2^2\|u_t\|_2^2 + \|Au\|_2^2\|\tht_0\|_2^2
\end{equation}
 Recall \eqref{e:u.t.avg}, Theorem~\ref{t:l2} and Theorem~\ref{t:v.l2} to
 conclude that $\int_t^{t+r}\|u_t(\tau)\|_2^2\ d\tau$, $\|\nbl u(t)\|_2$ and
 $\int_t^{t+r}\|Au(\tau)\|_2^2\ d\tau$ are all uniformly bounded. Moreover,
 there is a {\em bounded} absorbing set for all the three quantities.
 Therefore, we can apply Lemma~\ref{l:ug} to \eqref{e:u.t.2} to get the
 existence of a bounded set in $\R_+$ which absorbs $\|u_t(t)\|_2$. Moreover,
 with the help of the local regularity in {\bf Step 1}, we get the uniform
 boundedness of $\|u_t(t)\|_2$ for all $t\in(0,+\infty)$.
 Then, integrate \eqref{e:u.t.2} with respect to $t$ and use the fact
 $u_t\in L^\infty(0,\infty; H)$ and Therorem~\ref{t:v.l2} to get
 $u_t\in L^2_{loc}(0, \infty; V)$.

 {\bf Step 3}. Global regularity $u\in L^\infty(0, \infty; D(A))$. 
 
 By \eqref{e:u.v}, we have
\begin{equation*}
\begin{split}
 \|Au(t)\|_2^2 &= \lra{{\field{P}}\tht e_2 -u_t-B(u,u)}{Au}\\
&\ls (\|\tht\|_2 + \|u_t\|_2)\|Au\|_2
   + C\|u\|_2^\hf\|A^\hf u\|_2\|Au\|_2^{\frac{3}{2}}\\
&\ls C(\|\tht\|_2^2 + \|u_t\|_2^2 + \|u\|_2^2\|A^\hf u\|_2^4)
 +\hf \|Au(t)\|_2^2.
\end{split}
\end{equation*}
 Thus,
\begin{equation*}
\|Au(t)\|_2^2 \lsp \|\tht_0\|_2^2 + \|u_t\|_2^2 + \|u\|_2^2\|A^\hf u\|_2^4.
\end{equation*}
 Therefore, $\|Au(t)\|_2$ is also uniformly bounded for $t\in[0,\infty)$ and
 there exists a bounded absorbing set for $\|Au(t)\|_2$ in $\R_+$ as well.

\end{prf}

\subsection{Global Single Exponential Estimate of $\|\nbl\tht\|_2$}

 Now, we prove global boundedness of $\|\nbl\tht\|_2$. Our approach is still
 quite different from previous methods used in \cite{Hu;Kukavica;Ziane:2013}.
 Notice that, {\em even if assuming that the uniform boundedness of
 $\|Au(t)\|_2$ were already known}, following previous approache would still
 results in {\em double exponential growth} of the upper bound of
 $\|\nbl\tht(t)\|_2$ with respect to $t$.

\begin{thrm}
\label{t:t.h1}
 Suppose $(u_0,\tht_0)\in(D(A), H^1)$. Then, for any $T>0$,  we have
 $\tht \in L^\infty(0, T; H^1)$. Moreover, there is a constant $C>0$
 depending only on $t_1(>0)$, $u_0$ and $\tht_0$, such that 
\begin{equation}
\label{e:t.h1.bd}
 \|\nbl\tht(t)\| \ls C e^{Ct^2},\quad \forall t\gs t_1>0.
\end{equation}
 Moreover, there exists a $C>0$ depending only $\tht_0$, such that
\begin{equation}
\label{e:t.h1.asmp}
 \lim_{t\goto\infty}\|\nbl\tht(t)\|e^{-Ct^2} \ls C.
\end{equation}
\end{thrm}

\begin{prf}

 By the local regularity result Proposition~\ref{p:h2.h1.loc}, we only need
 to prove global boundedness of $\tht\in L^\infty_{\text{loc}}(t_1, \infty; H^1)$,
 or simply $\nbl\tht\in L^\infty_{\text{loc}}(t_1, \infty; L^2)$, and it is
 enough to prove \eqref{e:t.h1.bd}.

 Since $\|Au(t)\|_2$ is uniformly bounded for $t\in[0,+\infty)$ by
 Theorem~\ref{t:u.h2}, we can use Brezis-Gallouet inequality (see
 \cite{Brezis;Gallouet:1980}) to obtain
\begin{equation}
\label{e:B.G}
 \|\nbl u(t)\|_\infty \ls C\lf[1+\log^\hf(1+\|A^\frac{3}{2}u(t)\|_2)\rt],
\end{equation}
 where $C$ is {\em independent} of $t$ and $u$ and {\em depends only on}
 $\|\tht_0\|_2$. Therefore, by \eqref{e:t.nbl}, we get
\[ \dert{\|\nbl\tht\|_2^2} \ls
  C\lf[1+\log^\hf(1+\|A^\frac{3}{2}u(t)\|_2)\rt]\|\nbl\tht\|_2^2,\]
 that is
\[  \|\nbl\tht(t)\|_2^2 \ls \|\nbl\tht_0\|_2^2
\exp\lf\{C\int_0^t\lf[1+\log^\hf(1+\|A^\frac{3}{2}u(\tau)\|_2)\rt] d\tau\rt\}.\]
 Due to concavity of the function $f(x) := \log^\hf(1+x)$ for $x\in[0,+\infty)$,
 Jensen's inequality implies
\[ \frac{1}{t}\int_0^t\lf[1+\log^\hf(1+\|A^\frac{3}{2}u(\tau)\|_2)\rt]\ d\tau
 \ls
 1 + \log^\hf\lf(1+ \frac{1}{t}\int_0^t\|A^\frac{3}{2}u(\tau)\|_2\ d\tau\rt).\]
 Hence, the above two inequalities yield
\begin{equation*}
\begin{split}
\|\nbl\tht(t)\|_2^2 &\ls \|\nbl\tht_0\|_2^2
\exp\lf\{Ct + Ct\log^\hf\lf[1+ \frac{1}{t}\int_0^t\|A^\frac{3}{2}u(\tau)\|_2\
 d\tau\rt]\rt\}\\
&\ls \|\nbl\tht_0\|_2^2
 \lf[1+\frac{1}{t}\int_0^t\|A^\frac{3}{2}u(\tau)\|_2\ d\tau \rt] e^{C(t+t^2)}.
\end{split}.
\end{equation*}
 Notice that in the last step of the above inequality, we have used
 Cauchy-Schwartz inequality.

 Let $t\gs t_1>0$. Then, we obtain form the above estimate
\begin{equation*}
\begin{split}
\|\nbl\tht(t)\|_2^2 &\ls \|\nbl\tht_0\|_2^2 e^{C(t+t^2)} + 
\|\nbl\tht_0\|_2^2\frac{e^{C(t+t^2)}}{t}\int_0^t\|A^\frac{3}{2}u(\tau)\|_2\
  d\tau\\
&\ls \|\nbl\tht_0\|_2^2 e^{C(t+t^2)}
  + \|\nbl\tht_0\|_2^2\frac{e^{C(t+t^2)}}{\sqrt{t}}
 \lf(\int_0^t\|A^\frac{3}{2}u(\tau)\|_2^2\ d\tau\rt)^\hf\\
&\ls \lf[\|\nbl\tht_0\|_2^2 + \frac{\|\nbl\tht_0\|_2^4}{t_1}\rt]e^{C(t+t^2)}
 + \int_0^t\|A^\frac{3}{2}u(\tau)\|_2^2\ d\tau.
\end{split}
\end{equation*}
 Integrating \eqref{e:u.h2}, we have
\[ \int_0^t\|A^\frac{3}{2}u(\tau)\|_2^2\ d\tau
 \ls \|Au_0\|_2^2 + C\int_0^t\lf[ \|\tht_0\|_2^2 + \|u\|_2\|Au\|_2^3
 + \|\nbl\tht\|_2^2 \rt] d\tau. \]
Therefore, the above two inequalities give us, for $t\gs t_1$,
\begin{equation}
\label{e:t.nbl.1}
\|\nbl\tht(t)\|_2^2
 \ls\lf[\|\nbl\tht_0\|_2^2 + \frac{\|\nbl\tht_0\|_2^4}{t_1}\rt]e^{C(t+t^2)}
 + C(1+t) +C \int_0^t\|\nbl\tht(\tau)\|_2^2\ d\tau,
\end{equation}
 where we have used uniform boundedness of $u\in L^\infty(0,\infty,D(A))$ as
 given in Theorem~\ref{t:u.h2}. Denote
\[ y(t) := \int_0^t\|\nbl\tht(\tau)\|_2^2\ d\tau. \]
 Then, \eqref{e:t.nbl.1} becomes
\[ y'(\tau) = \|\nbl\tht(\tau)\|_2^2 \ls C(\tau+1)
 + C e^{C(\tau+\tau^2)} + C y(\tau), \]
 that is
\[ y'(\tau) - C y(\tau) \ls C(\tau+1) + Ce^{C(\tau+\tau^2)}. \]
 Integrating the above inequality for $\tau$ from $t_1$ to $t>t_1$, we get
\begin{equation}
\label{e:t.nbl.avg}
\begin{split}
 \int_0^t\|\nbl\tht(\tau)\|_2^2\ d\tau &= y(t)\\
&\ls e^{C(t-t_1)}\lf[ y(t_1) +C(t+1)(t-t_1)
  + C\int_{t_1}^t e^{C(\tau+\tau^2)}\ d\tau\rt].
\end{split}
\end{equation}
 Noticing that
\[ \int_{t_1}^te^{C(\tau+\tau^2)}\ d\tau
  \ls \int_{t_1}^te^{C(\tau+1)^2}\ d\tau \ls
  \frac{1}{Ct_1} \lf(e^{C(t+1)^2}-e^{C(t_1+1)^2}\rt) \ls Ce^{Ct^2}, \]
 and recalling local regularity result from Proposition~\ref{p:h2.h1.loc}:
\[ \|\nbl\tht(t)\|_2^2 \ls 2(\|Au_0\|_2^2 + \|u_0\|_2^2),\quad
  \forall t\in[0,t_1], \]
 we obtain from \eqref{e:t.nbl.avg} that, for $t\gs t_1>0$,
\begin{equation}
\label{e:t.nbl.avg.1}
  \int_0^t\|\nbl\tht(\tau)\|_2^2\ d\tau  \ls C e^{Ct^2},
\end{equation}
 where the $C$'s in the right-hand side of \eqref{e:t.nbl.avg.1} may depend only
 on $\|\tht_0\|_2$, $\|u_0\|$, $\|Au\|_0$ and $t_1$.  Now, \eqref{e:t.h1.bd}
 follows immediately from \eqref{e:t.nbl.1} and \eqref{e:t.nbl.avg.1}. It is
 easy to see from our proof that \eqref{e:t.h1.asmp} is also valid.
\end{prf}

\section{Global regularity in $V\times H^1$}
\label{s:h1.h1} 

 In this section, we prove global regularity of the solution $(u,\tht)\in
 V\times H^1$. As can be seen from Section~\ref{s:h2.h1}, the Brezis-Gallouet
 inequality \eqref{e:B.G} is one of the key steps in obtaining estimate on
 $\|\nbl\tht(t)\|_2$. However, to use \eqref{e:B.G}, one has to deal with
 $\|A^{\frac{3}{2}}u\|_2$. Therefore, it is natrual to introduce the estimate
 \eqref{e:u.h2}. This is the main reason for studying local and global
 regularity of $(u,\tht)\in D(A)\times H^1$. Notice that, it is possible to use
 a less demanding variant of Brezis-Gallouet inequality involving only
 $\|A^{1+\veps}u\|_2$, with $\veps>0$, which however would be quite complicated.
 If we now study local and global regularity of $(u,\tht)\in V\times H^1$, the
 estimate \eqref{e:u.h2} will not be available and more importantly the
 Brezis-Gallouet type inequality \eqref{e:B.G} (and its variants) can not help.
 Therefore, {\em essential new difficulties appear here}. It is quite remarkable
 that our following analysis using the idea of spectral decompostion can help to
 resolve the problem rather conveniently. Our next main result of this article
 is presented the following:

\begin{thrm}
\label{t:h1}
 Suppose $(u_0,\tht_0)\in V\times H^1$ and $(u,\tht)$ is the {\em (unique)} weak
 solution to the equations \eqref{e:u}-\eqref{e:t} with conditions \eqref{e:bc}
 and \eqref{e:ic}. Then,
\[ u\in L^\infty(0, \infty; V),\quad
   \tht\in L^\infty_{\text{loc}}(0,\infty; H^1). \]
 Moreover, for fixed $\tht_0$, there is a bounded absorbing set in $\R_+$ for
 $\|\nbl u(t)\|_2$ and $\|A u(t)\|_2$ as $t\goto\infty$, which is independent
 of $u_0$. Furthermore, there is a constant $C>0$ depending only on $t_1(>0)$,
 $u_0$ and $\tht_0$, such that 
\begin{equation}
\label{e:t.h1.bd.2}
 \|\nbl\tht(t)\| \ls C e^{Ct^2},\quad \forall t\gs t_1>0.
\end{equation}
\end{thrm}

\begin{prf}

{\bf Step 1}. We first prove local in time existence of the solution $(u,\tht)$
 in $V\times H^1$,  that is, there exists a $t_0=t_0(u_0,\tht_0)>0$, such that
 $(u,\tht)\in L^\infty(0, t_0; V\times H^1)$.

 Since $(u,\tht)\in L^\infty(0,\infty; V\times L^2)$, we only need to prove the
 local regularity of $\tht\in L^\infty(0, t_0; H^1)$ for some $t_0>0$. Due to
 existence and uniqueness of the weak solution, we only need to carry out some
 {\em a priori} estimates formally in the following proof. These {\em a priori}
 estimates can be rigorously justified with a standard approximation argument.

 By \eqref{e:t.nbl}, we have
\begin{equation}
\label{e:t.nbl.x}
 \|\nbl\tht(t)\|_2 \ls
 \|\nbl\tht_0\|_2 \exp\{\int_0^t\|\nbl u(\tau)\|_\infty\ d\tau\}
\end{equation}
 
 Apply $P_m$ to \eqref{e:u.v} and then take inner product with $P_mu$. We get
\[ \hf\dt{\|P_m u(t)\|_2^2} + \lra{P_mA u}{P_m u}
 = \lra{P_m\lf[\Pf(\tht e_2) -B(u,u)\rt]}{P_m u}. \]
 Noticing that
\[ \lra{P_mA u}{P_mu} = \sum_{2^{m-1}\ls j<2^m}|(u, w_j)|^2
 \gs 2^{m-1}\|P_m u\|_2^2,\quad \forall m\gs 1, \]
 we get
\[ \dt{\|P_m u(t)\|_2} + 2^{m-1}\|P_m u\|_2 \ls
 \|P_m\Pf(\tht e_2)\|_2 +\|P_mB(u,u)\|_2.\]
 Thus,
\begin{equation}
\label{e:u.pm}
\begin{split}
 2\|P_mu(t)\|_2+&\int_0^t 2^m\|P_m u(\tau)\|_2\ d\tau\\
 &\ls  2\|P_m u_0\|_2 + 2\int_0^t
  (\|P_m\Pf(\tht(\tau)e_2)\|_2+ \|P_m B(u(\tau), u(\tau))\|_2)\ d\tau.
\end{split}
\end{equation}
 Multiplying \eqref{e:u.pm} by $2^{\frac{m}{2}}$ and using \eqref{e:a.eqv} yields
\begin{equation}
\label{e:pm}
\begin{split}
\int_0^t 2^{m+\frac{m}{2}}\|P_m u(\tau)\|_2 d\tau
\lsp & \|A^\hf P_m u_0\|_2
  + \int_0^t\|A^\hf P_m \Pf(\tht(\tau)e_2)\|_2\ d\tau\\
& + \int_0^t\|A^\hf P_m B(u(\tau), u(\tau))\|_2\ d\tau.
\end{split}
\end{equation}
 Squaring both sides of \eqref{e:pm} and summing up with respect to $m$ yields
\begin{equation*}
\begin{split}
 \sum_{m=1}^\infty & \lf(\int_0^t 2^{m+\frac{m}{2}}\|P_mu\|_2\ d\tau\rt)^2\\
\lsp& \sum_{m=1}^\infty\|A^\hf P_m u_0\|_2^2
 +\sum_{m=1}^\infty \lf(\int_0^t\|A^\hf P_m\Pf(\tht(\tau)e_2)\|_2\ d\tau\rt)^2\\
 & +\sum_{m=1}^\infty
  \lf(\int_0^t\|A^\hf P_mB(u(\tau), u(\tau))\|_2\ d\tau\rt)^2\\
=&\|A^\hf u_0\|_2^2
 +\sum_{m=1}^\infty \lf(\int_0^t\|A^\hf P_m\Pf(\tht(\tau)e_2)\|_2\ d\tau\rt)^2\\
 & +\sum_{m=1}^\infty
  \lf(\int_0^t\|A^\hf P_m B(u(\tau), u(\tau))\|_2\ d\tau\rt)^2.
\end{split}
\end{equation*}
 Then, take square root of the above inequality to get
\begin{equation}
\label{e:sqr.sm.pm}
\begin{split}
&
\lf[\sum_{m=1}^\infty\lf(\int_0^t 2^{m+\frac{m}{2}}\|P_mu\|_2\ d\tau\rt)^2\rt]^\hf
\\ 
\lsp & \|A^\hf u_0\|_2 + \lf[\sum_{m=1}^\infty
 \lf(\int_0^t\|A^\hf P_m\Pf(\tht(\tau)e_2)\|_2\ d\tau\rt)^2 \rt]^\hf \\
 & + \lf[\sum_{m=1}^\infty
 \lf(\int_0^t\|A^\hf P_m B(u(\tau), u(\tau))\|_2\ d\tau\rt)^2\rt]^\hf \\
\ls &\|A^\hf u_0\|_2 + \int_0^t \lf(\sum_{m=1}^\infty
 \|A^\hf P_m\Pf(\tht(\tau)e_2)\|_2^2 \rt)^\hf\ d\tau \\
 & + \int_0^t
 \lf(\sum_{m=1}^\infty\|A^\hf P_m B(u(\tau), u(\tau))\|_2^2 \rt)^\hf\ d\tau\\
=& \|A^\hf u_0\|_2  + \int_0^t \|A^\hf\Pf(\tht(\tau)e_2)\|_2\ d\tau
 + \int_0^t \|A^\hf B(u(\tau), u(\tau))\|_2\ d\tau,
\end{split}
\end{equation}
 where Minkowski's inequality is used twice for the second inequality of
 \eqref{e:sqr.sm.pm}.

 By Agmon's inequality, \eqref{e:a.eqv} and \eqref{e:a.hf}, it is easy to see
 that
\begin{equation}
\label{e:dlt.pm}
\begin{split}
 \|\nbl P_m u(\tau)\|_\infty
\lsp& \|\nbl P_m u(\tau)\|_2^\hf\|A\nbl P_m u(\tau)\|_2^\hf\\
\lsp& 2^{\frac{m}{2}}\|\nbl P_mu(\tau)\|_2
= 2^{\frac{m}{2}}\|A^\hf P_mu(\tau)\|_2 \lsp 2^m\|P_mu(\tau)\|_2.
\end{split}
\end{equation}
 Therefore, by \eqref{e:dlt.pm} and \eqref{e:sqr.sm.pm}, we obtain
 \begin{equation}
\begin{split}
\label{e:u.nbl.inf}
 \int_0^t\|\nbl u(\tau)\|_\infty\ d\tau
\ls& \int_0^t\sum_{m=1}^\infty\|\nbl P_m u(\tau)\|_\infty\ d\tau\\
=& \sum_{m=1}^\infty 2^{-\frac{m}{2}}
 \int_0^t 2^{m+\frac{m}{2}}\|P_m u(\tau)\|_2\ d\tau\\
\ls & \lf[\sum_{m=1}^\infty
   \lf( \int_0^t 2^{m+\frac{m}{2}}\|P_mu\|_2\ d\tau \rt)^2 \rt]^\hf\\
\lsp& \|\nbl u_0\|_2 + \int_0^t \|\nbl\Pf(\tht(\tau)e_2)\|_2\ d\tau
 + \int_0^t \| \nbl B(u(\tau), u(\tau))\|_2\ d\tau.
\end{split}
\end{equation}
 The last term of \eqref{e:u.nbl.inf} can be estimated using \eqref{e:b} as
\begin{equation}
\begin{split}
\label{e:b.nbl}
 \int_0^t\|\nbl B(u(\tau),u(\tau))\|_2\ d\tau
&\lsp \int_0^t \|u(\tau)\|_2^\hf \|Au(\tau)\|_2^{\frac{3}{2}}\ d\tau \\
&\lsp \lf(\int_0^t \|u(\tau)\|_2^2\ d\tau\rt)^{\frac{1}{4}}
\lf(\int_0^t\|Au(\tau)\|_2^2\ d\tau \rt)^{\frac{3}{4}}.
\end{split}
\end{equation}
 Combining \eqref{e:u.nbl.inf} and \eqref{e:b.nbl}, we get
\begin{equation}
\label{e:u.nbl.inf-1}
\begin{split}
  \int_0^t\|\nbl u(\tau)\|_\infty\ d\tau
\lsp& \|\nbl u_0\|_2 + t + \int_0^t\|\nbl\tht(\tau)\|_2\ d\tau\\
\lsp& \|\nbl u_0\|_2  + \int_0^t\lf(1+\|\nbl\tht(\tau)\|_2\rt) d\tau,
\end{split}
\end{equation}
 where we have used Theorem~\ref{t:l2}, Theorem~\ref{t:v.l2} and
 \eqref{e:u.nbl.l2.1}. By \eqref{e:t.nbl.x} and \eqref{e:u.nbl.inf-1}, we obtain
\begin{equation}
\label{e:t.nbl.x.2}
 \|\nbl\tht(t)\|_2
 \lsp \exp\lf\{c_0\int_0^t\lf(1+\|\nbl\tht(t)\|_2\rt) d\tau\rt\},
\end{equation}
 where $c_0>0$ is a constant depending on initial data.
 Denote:
\[ y(t):=c_0\int_0^t\lf(1+\|\nbl\tht(t)\|_2\rt) d\tau.\]
 By \eqref{e:t.nbl.x.2}, we have
\[ y'(t) \ls C_1 e^{y(t)}, \]
 where $C_1>0$ is a constant depending on $u_0$ and $\tht_0$. Thus, there
 exists a $t_0>0$, such that
\[ c_0\int_0^t\lf(1+\|\nbl\tht(t)\|_2\rt) d\tau = y(t)
 \ls \ln\lf(\frac{1}{1-C_1t}\rt) < +\infty ,\quad \text{for } t\in[0,t_0]. \]
 Therefore, by \eqref{e:t.nbl.x.2},
\[ \|\nbl\tht(t)\|_2 \ls \frac{C_2}{1-C_1 t}<+\infty,\quad t\in[0, t_0], \]
 where $C_2>0$ is also a constant depending on $u_0$ and $\tht_0$. This finishes
 the proof of $\tht\in L^\infty(0, t_0; H^1)$.

 {\bf Step 2}. We prove global in time existence of $\tht$ in $H^1$.

 Suppose that $(u_0,\tht_0)\in V\times H^1$. Then, by {\bf Step 1}, there is
 $t_0>0$, such that $(u,\tht)\in L^\infty(0, t_0; V\times H^1)$.
 Due to Theorem~\ref{t:v.l2}, we can assume without losing generality that
 $u(t_0)\in D(A)$. Therefore, since $(u(t_0),\tht(t_0))\in D(A)\times H^1$, by
 the global existence result Theorem~\ref{t:u.h2} and Theorem~\ref{t:t.h1},
 we have $u\in L^\infty(t_0,+\infty; D(A))$ and
 $\tht\in L^\infty_{\text{loc}}(t_0,\infty; H^1)$, thus proving the
 global regularity $(u,\tht)\in L^\infty(0,T; V\times H^1)$ for any $T>0$.
 The rest part of the theorem follows from Theorem~\ref{t:v.l2},
 Theorem~\ref{t:u.h2} and Theorem~\ref{t:t.h1} as well.

\end{prf}

\begin{rmrk}
 Theorem~\ref{t:h1} is obviously also valid for $\Omg=\T^2$. By Littlewood-Paley
 frequence decomposition, the method of our proof of global regularity in
 $V\times H^1$ can be extended immediately to the case of $\Omg=\R^2$.
\end{rmrk}

\section{Global Regularity in $D(A)\times H^2$}
\label{s:h2.h2} 

 The following theorem is our last main result for global regularity, which
 we include here for completeness. Since we will not use Brezis-Gallouet type
 inequalities in our proof of Theorem~\ref{t:h2}, the proof is a lot simpler
 than otherwise.

\begin{thrm}
\label{t:h2}
 Suppose $(u_0,\tht_0)\in D(A)\times H^2$ and $(u,\tht)$ is the {\em (unique)}
 weak solution to the equations \eqref{e:u}-\eqref{e:t} with conditions
 \eqref{e:bc} and \eqref{e:ic}. Then,
\[ u\in L^\infty(0, \infty; D(A)),\quad \text{and} \quad 
   \tht\in L^\infty([0,T]; H^2), \ \forall T>0. \]
 Moreover, for fixed $\tht_0$, there is a bounded absorbing set in $\R_+$ for
 $\|\nbl u(t)\|_2$ and $\|A u(t)\|_2$ as $t\goto\infty$, which is independent
 of $u_0$.
\end{thrm}

\begin{prf}

 It is easy to see that we just need to prove that for any $T>0$
 \[ \tht\in L^\infty(0, T; H^2). \]
 All the other statements of the theorem are already confirmed by the previous
 theorems.
 
 Applying $\prt_x^2$ to \eqref{e:t} and then taking inner product of the
 derived equation with $\prt_x^2\tht$, we obtain
\begin{equation*}
\begin{split}
 \dert{\|\prt_x^2\tht\|_2^2}
&=-2\lra{(\prt_x^2u)\cdot\nbl\tht
  +2(\prt_x u)\cdot\nbl(\prt_x\tht)}{\prt_x^2\tht} \\
& \lsp \|\nbl\tht\|_4\|A u\|_4\|\nbl^2\tht\|
 + \|\nbl u\|_\infty\|\nbl^2\tht\|_2^2\\
& \lsp \|\nbl\tht\|_2^\hf\|\nbl^2\tht\|_2^{\frac{3}{2}}
  \|A u\|_2^\hf\|A^{\frac{3}{2}} u\|_2^\hf +\|\nbl u\|_\infty\|\nbl^2\tht\|_2^2\\
&\lsp \|\nbl\tht\|_2^2\|Au\|_2^2\|A^{\frac{3}{2}}u\|_2^2
 + (1+\|\nbl u\|_\infty)\|\nbl^2\tht\|_2^2.
\end{split}
\end{equation*}
 Similar estimates can be obtained for $\prt_{xy}^2\tht$ and $\prt_y^2\tht$.
 Therefore, we have
\begin{equation*}
 \dert{\|\nbl^2\tht\|_2^2}
 \lsp \|\nbl\tht\|_2^2\|Au\|_2^2\|A^{\frac{3}{2}}u\|_2^2
 + (1+\|\nbl u\|_\infty)\|\nbl^2\tht\|_2^2.
\end{equation*}
 Applying Gronwall lemma to the above inequality, we obtain
\begin{equation}
\label{e:t.h2}
\begin{split}
 \|\nbl^2\tht(t)\|_2^2
\lsp& \|\nbl^2\tht_0\|_2^2 e^{\int_0^t(1+\|\nbl u(\tau)\|_\infty)\ d\tau}\\
 &+ \int_0^t
 \|\nbl\tht(\tau)\|_2^2\|Au(\tau)\|_2^2\|A^{\frac{3}{2}}u(\tau)\|_2^2
 e^{\int_\tau^t(1+\|\nbl u(s)\|_\infty)\ ds} \ d\tau.
\end{split}
\end{equation}
 By \eqref{e:u.nbl.inf} and \eqref{e:b.nbl}, we have
\begin{equation}
\label{e:nbl.u.inf}
\begin{split}
 \int_0^t \|\nbl u(\tau)\|_\infty\ d\tau
\lsp& \|\nbl u_0\|_2 + \int_0^t \|\nbl\tht(\tau)\|_2\ d\tau\\
 &+ \lf(\int_0^t\|u(\tau)\|_2^2\ d\tau\rt)^{\frac{1}{4}}
   \lf(\int_0^t\|Au(\tau)\|_2^2\ d\tau\rt)^{\frac{3}{4}}
\end{split}
\end{equation}
 {\bf Case 1}. Suppose that $t\in[0, t_1]$.
 Then, by \eqref{e:nbl.u.inf}
\[ \int_0^t \|\nbl u(\tau)\|_\infty\ d\tau \ls
  \int_0^{t_0} \|\nbl u(\tau)\|_\infty\ d\tau
 \ls C < \infty,\]
 where we have applied Theorem~\ref{t:l2}, Theorem~\ref{t:v.l2} (or
 Theorem~\ref{t:u.h2}) and Theorem~\ref{t:t.h1} to \eqref{e:nbl.u.inf}
 with $t$ replaced by $t_0$. Then, applying these theorems and the above
 inequality to \eqref{e:t.h2}, we have
\[ \|\nbl^2\tht(t)\|_2^2 \lsp \|\nbl^2\tht_0\|_2^2 + 
 \int_0^{t_0}\|A^{\frac{3}{2}}u(\tau)\|_2^2\ d\tau,\quad \forall t\in[0, t_1]. \]
 Then, by \eqref{e:u.h2}, we get that for $t\in[0, t_1]$
\[ \|\nbl^2\tht(t)\|_2^2 \lsp \|\nbl^2\tht_0\|_2^2 + \|\tht_0\|_2^2 t
 + \int_0^t(\|\nbl\tht(\tau)\|_2^2+\|u(\tau)\|_2^2\|Au(\tau)\|_2^2)\ d\tau, \]
 from which we derive regularity of local boundedness:
\[ \tht \in L^\infty(0, t_1, H^2). \]

 {\bf Case 2}. Suppose $t\in [t_1, T]$. This can be treated the same way
 as {\bf Case 1} and we can obtain a upper bound on $\|\nbl^2\tht(t)\|_2^2$
 for $t\in [t_1, T]$.

\end{prf}

\begin{rmrk}
 Since \eqref{e:u.nbl.inf} and \eqref{e:b.nbl} are all still valid for
 $\Omg=\R^2$. Global regularity of $(u,\tht)\in D(A)\times H^2$ is also still
 valid when $\Omg=\R^2$.
\end{rmrk}


\begin{thebibliography}{99}

\bibitem{Brezis;Gallouet:1980} H. Brezis and T. Gallouet,
 {\em Nonlinear Schr¡\"odinger evolution equations},
 Nonlinear Anal. Theory, Methods Appl. {\bf 4} (1980), No.4, 677-681.

\bibitem{Chae:2006} D. Chae,
 {\em Global regularity for the 2D Boussinesq equations with partial viscosity terms}, Adv. Math. {\bf 203}, (2006), No. 2, 497-513.

\bibitem{Constantin;Foias:1988} P. Constantin and C. Foias,
 Navier-Stokes Equations, The University of Chicago Press, 1988.

\bibitem{Danchin;Paicu:2008} R. Danchin and M. Paicu,
 {\em Le th\'eor\`eme de Leray et le th\'eor\`eme de Fujita-Kato pour le syst\`eme de Boussinesq partiellement visqueux}
, Bull. Soc. Math. France {\bf 136} (2008), no. 2, 261-309.

\bibitem{Danchin;Paicu:2011} R. Danchin and M. Paicu,
 {\em Global existence results for the anisotropic Boussinesq system in dimension two}, Math. Models Meth. Appl. Sci. {\bf 21} (2011), No.3, 421-457.

\bibitem{He:2012} L. He,
 {\em Smoothing estimates of 2d incompressible Navier-Stokes equations in bounded domains with applications},
 J. Funct. Anal. {\bf 262} (2012) No. 7, 3430-3464.

\bibitem{Hu;Kukavica;Ziane:2013} W. Hu, I. Kukavica and M. Ziane,
 {\em On the regularity for the Boussinesq equations in a bounded domain},
 J. Math. Phys., {\bf 54} (2013), no. 8, 081507, 10pp.

\bibitem{Hou;Li:2005} T. Y. Hou and C. Li,
 {\em Global well-posedness of the viscous Boussinesq equations},
 Discrete Contin. Dyn. Syst. {\bf 12} (2005), No. 1, 1-12.

\bibitem{Lai;Pan;Zhao:2011} M. J. Lai, R. Pan, and K. Zhao,
 {\em Initial boundary value problem for two-dimensional viscous Boussinesq equations}, Arch. Ration. Mech. Anal. {\bf 199}(3), 739-760 (2011).

\bibitem{Larios;Lunasin;Titi:2013} A. Larios, E. Lunasin, and E. S. Titi,
 {\em Global well-posedness for the 2D Boussinesq system with anisotropic viscosity and without heat diffusion}, 
 J. Differ. Equations {\bf 255} (2013), No. 9, 2636-2654.

\bibitem{LiYG:2004} Y. C. Li,
 {\em Global regularity for the viscous Boussinesq equations},
 Math. Methods Appl. Sci. {\bf 27} (2004), no. 3, 363-369. 

\bibitem{Temam:1977} R. Temam,
 {\em Navier-Stokes equations. Theory and numerical analysis},
 Studies in Mathematics and its Applications, Vol. 2. North-Holland, 1977.
 Reprinted with corrections by AMS, 2001.

\bibitem{Temam:1988} R. Temam,
 {\it Infinite Dimensional Dynamical Systems in Mechanics and Physics}
 Spring-Verlag, New York, Applied Mathematical Sciences Series, vol. 68, 1988.
 Second augmented edition, 1997.

\end{thebibliography}
\end{document}